\documentclass[12pt]{article}
\newcommand{\ncr}[2]{\mbox{$\left(\begin{array}{c}#1\\#2\end{array}\right)$}}
\newcommand{\bc}{\begin{center}}
\newcommand{\ec}{\end{center}}
\newcommand{\be}{\begin{equation}}
\newcommand{\ee}{\end{equation}}
\newcommand{\bea}{\begin{eqnarray}}
\newcommand{\eea}{\end{eqnarray}}
\newcommand{\bfl}{\begin{flushleft}}
\newcommand{\efl}{\end{flushleft}}
\newcommand{\bdm}{\begin{displaymath}}
\newcommand{\edm}{\end{displaymath}}
\newcommand{\ba}{\begin{array}}
\newcommand{\ea}{\end{array}}
\newcommand{\bd}{\begin{description}}
\newcommand{\ed}{\end{description}}
\newcommand{\ben}{\begin{enumerate}}
\newcommand{\een}{\end{enumerate}}
\newcommand{\beas}{\begin{eqnarray*}}
\newcommand{\eeas}{\end{eqnarray*}}
\newcommand{\bb}{\begin{thebibliography}{99}}
\newcommand{\eb}{\end{thebibliography}}
\newcommand{\bs}{\begin{sloppypar}}
\newcommand{\es}{\end{sloppypar}}

\newtheorem{claim}{Claim}[section]
\newtheorem{conjecture}[claim]{Conjecture}  
\newtheorem{definition}[claim]{Definition}
\newtheorem{lemma}[claim]{Lemma}
\newtheorem{theorem}[claim]{Theorem}
\newtheorem{corollary}[claim]{Corollary}
\newtheorem{proposition}[claim]{Proposition} 
\newcommand{\bclaim}{\begin{claim}}
\newcommand{\eclaim}{\end{claim}}
\newcommand{\blem}{\begin{lemma}}
\newcommand{\elem}{\end{lemma}}
\newcommand{\bthm}{\begin{theorem}}
\newcommand{\ethm}{\end{theorem}}
\newcommand{\bprop}{\begin{proposition}}
\newcommand{\eprop}{\end{proposition}} 
\newcommand{\bcor}{\begin{corollary}}
\newcommand{\ecor}{\end{corollary}}
\newcommand{\bcon}{\begin{conjecture}}
\newcommand{\econ}{\end{conjecture}}

\title{A reconstruction problem related to balance equations-II: the 
general case }
\author{Bhalchandra D. Thatte \\
\small Allan Wilson Centre for Molecular Ecology and Evolution, \\[-0.8ex]
\small and Institute of Fundamental Sciences, \\[-0.8ex]
\small Massey University, Palmerston North, New Zealand \\[-0.8ex]
\small \texttt{b.thatte@massey.ac.nz}}

\date{\small Discrete Mathematics 194, no. 1-3(1999) 281-284. \\
\small Mathematics Subject Classifications: 05C60}
\begin{document}
\maketitle
\begin{abstract}
A modified $k$-deck of a graph $G$ is obtained by removing $k$
edges of $G$ in all possible ways, and adding $k$
(not necessarily new) edges in all possible ways. Krasikov
and Roditty asked if it was possible to construct the usual
$k$-edge deck of a graph from its modified $k$-deck. Earlier
I solved this problem for the case when $k=1$. In this paper, 
the problem is completely solved for arbitrary $k$.
The proof makes use of the $k$-edge version of Lov\'asz's result and 
the eigenvalues of certain matrix related to the Johnson graph.

This version differs from the published version.
Lemma 2.3 in the published version had a typo in one
equation. Also, a long manipulation of some combinatorial expressions
was skipped in the original proof of Lemma 2.3,
which made it difficult to follow the proof.
Here a clearer proof is given.

\end{abstract}

\section{Introduction}
\label{intro}
The graphs considered in this paper are simple and undirected, and 
are assumed to have $n$ vertices. The complement of $G$ is denoted by $G^c$. 
Let $N = \ncr{n}{2}$.
Let $U_m$ denote the collection of all unlabelled $n$-vertex, $m$-edge
graphs. We define three matrices $\Delta_i$, $D_i$ and $d_i$ as follows.
The rows and columns of $\Delta_i$ and $D_i$ are indexed by the members of
$U_m$.  The $kl$-th entry of $\Delta_i$ is the number of graphs isomorphic
to $G_k$ that can be obtained by removing $i$ edges from $G_l$ and then
adding $i$ edges. Here the added edges need not be different from the
removed edges. The entries of $D_i$ are similarly defined with an
additional condition that the removed set of edges and the added set of
edges be disjoint. The rows of $d_i$ are indexed by $F_k \in U_{m-i}$, and
its columns are indexed by $G_l \in U_m$.  The $kl$-th entry of $d_i$ is
the number of $i$-edge deleted subgraphs of $G_l$ that are isomorphic to
$F_k$.  A set (or a multiset) $P$ of $m-i$-edge graphs is
denoted by its characteristic vector $X_P$ of length equal to $|U_{m-i}|$.
The characteristic vector of a singleton set $\{G\}$ is denoted by simply
$X_G$. This has only one entry equal to $1$ and other entries equal to
$0$. Thus, in our notation, the vector $d_kX_G$ represents the $k$-edge
deck of $G$, (denoted by $k-ED(G)$), and the vector $\Delta_kX_G$
represents the modified $k$-deck of $G$, i.e., the collection of graphs
obtained from $G$ by removing $k$ edges and then adding $k$ (not
necessarily new) edges.

Krasikov and Roditty first introduced modified decks for the 
purpose of proving the reconstruction result of M\"uller. 
They asked if  the $k$-edge deck of a graph could be constructed 
from  its modified $k$-deck.
In our notation, it is equivalent to asking if the vector $d_kX_G$
could be computed given the vector $\Delta_kX_G$.
In [T], this problem was solved for the case when $k=1$. 
Two proofs of 
this were offered there. In one proof, it was demonstrated that 
$\Delta_iX_G$ could be computed for $i> 1$ given $\Delta_1X_G$. 
The rest of the 
proof was based on the fact that Lov\'asz's edge reconstruction 
result in $k=1$ case could be proved directly from modified decks, i.e., 
without knowing the $1$-edge deck.
In the second proof, which was based on the eigen values of Johnson graph,
it was shown that Lov\'asz's result could be proved 
directly from $\Delta_1X_G$, thus avoiding the explicit construction 
of $\Delta_iX_G, \, i> 1$ in terms of $\Delta_1X_G$. 

The proof for the general case presented here does involve 
construction of $\Delta_iX_G$ in terms of $\Delta_kX_G$, for $i\geq k$.
But rest of the proof makes use of eigenvalues of Johnson graph.

\section{Reconstructing $d_kX_G$ from $\Delta_kX_G$ }
\label{mainresult}

In the following, we assume that for two graphs $G$ and $H$, 
we are given that $\Delta_kX_G = \Delta_kX_H$. We write $X=X_G -X_H$,
therefore, $\Delta_kX=0$. 
We first state two identities without proof. 
The first one - Lemma~\ref{DelD} - is 
equivalent to Lemma 3.1 in [KR], and the second one - 
Lemma~\ref{recursion} - is Theorem 2.2 from [T].  

\blem
\label{DelD}
$\Delta_s = \sum_{i=0}^s \ncr{m-i}{s-i}D_i $
\elem

\blem 
\label{recursion}
$D_1D_i=(m-i+1)(N-m-i+1)D_{i-1}
+ i(N-2i)D_i + (i+1)^2D_{i+1}$.
\elem  

\blem 
\label{DelDel}
\beas
\Delta_{i+1} = \frac{1}{(i+1)^2}
\left\{i(2m-N-i-1)\Delta_0 + \Delta_1\right\}\Delta_i
\eeas 
\elem 

\noindent {\bf Proof}
From Lemma ~\ref{recursion} we write
\bdm
(i+1)^2D_{i+1} = D_1D_i - (m-i+1)(N-m-i+1)D_{i-1} - i(N-2i)D_i
\edm
Substituting for $D_{i+1}$ and $D_i$ from Lemma~\ref{DelD},
we have
\begin{eqnarray}
\lefteqn{(i+1)^2\left(\Delta_{i+1}-
\sum_{j=0}^{i}\ncr{m-j}{i+1-j}D_j\right)}\nonumber \\[5mm]
&=& D_1\left( \Delta_i - \sum_{j=0}^{i-1}\ncr{m-j}{i-j}D_j\right) 
-(m-i+1)(N-m-i+1)D_{i-1}  \nonumber \\[5mm]
&-& i(N-2i)\left(\Delta_i - \sum_{j=0}^{i-1}\ncr{m-j}{i-j}D_j\right)
\nonumber
\end{eqnarray}
In the first term on the RHS, we substitute
$D_1\Delta_i = (\Delta_1 - m\Delta_0)\Delta_i$,
$D_1D_j;\,j>0$ from Lemma~\ref{recursion}, and $D_1D_0=D_1$.
Therefore,
\newpage
\begin{eqnarray}
\lefteqn{(i+1)^2\Delta_{i+1}}\nonumber \\[5mm]
& = & (\Delta_1 - m\Delta_0)\Delta_i -\ncr{m}{i}D_1 \nonumber \\[5mm]
& - & \sum_{j=1}^{i-1} \ncr{m-j}{i-j}(m-j+1)(N-m-j+1)D_{j-1} \nonumber \\[5mm]
& - & \sum_{j=1}^{i-1} \ncr{m-j}{i-j}j(N-2j)D_j  \nonumber \\[5mm]
& - & \sum_{j=1}^{i-1} \ncr{m-j}{i-j}(j+1)^2D_{j+1} \nonumber \\[5mm]
&-&(m-i+1)(N-m-i+1)D_{i-1}  \nonumber \\[5mm]
&-& i(N-2i)\left(\Delta_i - \sum_{j=0}^{i-1}\ncr{m-j}{i-j}D_j\right)
+(i+1)^2\sum_{j=0}^{i}\ncr{m-j}{i+1-j}D_j \nonumber
\end{eqnarray}
Two terms on the RHS contribute to $D_i$ -
the summation in the fourth line on the RHS, for $j=i-1$, and the last
summation in the last line on the RHS, for $j=i$.
Both these $D_i$ terms are replaced by
$\Delta_i - \sum_{j=0}^{i-1}\ncr{m-j}{i-j}D_j$.
This leaves only terms containing $D_j;\,j \leq i-1$.
One can then verify that, after simplification of the RHS, all
terms containing $D_j;\,j\leq i-1$ cancel out, and we get
\bdm
(i+1)^2\Delta_{i+1} = \left(i(2m-N-i-1)\Delta_0 + \Delta_1\right)\Delta_i
\edm
This completes the proof.

\bcor
\label{cor-DelDel}
If $\Delta_kX=0$ then $\Delta_iX=0$ for 
all $i \geq k$.
\ecor  

The following lemma is the $k$-edge version of 
Lov\'asz's result. This may be found in [GKR], but we only note here 
that the bound in the following result doesn't depend upon 
the number of graphs in the collection $P$.  
\blem
\label{lov}
Let $2p-k+1 > N$, and let $P$ and $Q$ be collections of 
$p$-edge graphs such that $d_kX_P=d_kX_Q$, then 
then  $X_P=X_Q$.   
\elem

Now we prove the main result of this section.
\bthm 
\label{Deld}
For collections $P$ and $Q$ of graphs, if $\Delta_kX_P=\Delta_kX_Q$ then 
$d_kX_P=d_kX_Q$. 
\ethm 

\noindent {\bf Proof} 
This is done by induction on $k$. The result was proved in [T] for 
$k=1$. Let the result be true
for $k \leq r-1$. 
Let $P' = \{F^c; F\in r-ED(P)\}$ and 
$Q' = \{F^c; F\in r-ED(Q) \}$. 
Here $r-ED(P)$ denotes the multiunion of $r$-edge decks 
of graphs in $P$. Note that  $\Delta_rX_P=\Delta_rX_Q$
is equivalent to  $d_rX_{P'}=d_rX_{Q'}$.
This follows from the fact that  
for any $F$, $A\in E(F)$ and $B$ disjoint with $E(F)-A$, 
$(F-A+B)^c = (F-A)^c -B$.
Now, if 
$2(N-m+r) -r +1 > N$, then $X_{P'}=X_{Q'}$, 
and $d_rX_P=d_rX_Q$. Therefore, we assume the contrary that $2m-r-1 \geq N$,
i.e., $2m-r+1 \geq N+2$. 

Now we demonstrate that either 
$\Delta_{r-1}X_G = \Delta_{r-1}X_H$ or $2m-r+1 \leq N+1$. 
We write, 
\beas
\Delta_r = \frac{1}{r^2}
\left\{(r-1)(2m-N-r)\Delta_0 + \Delta_1\right\}\Delta_{r-1}
\eeas 
We are interested in the invertibility of $(r-1)(2m-N-r)\Delta_0 + \Delta_1$.

\begin{definition}
Johnson graph is a simple graph whose vertex set is 
the family of $m$-sets
of an $N$-set. Two vertices $U$ and $V$ are adjacent if and only 
if $|U \cap V|=m-1$.
\end{definition}
 
Let $J$ be the adjacency matrix of the Johnson graph with parameters
$N=\ncr{n}{2}$ and $m$.  Let the square matrix $B$ be defined as follows.
The rows and columns of $B$ are indexed by all the labelled $m$-edge
graphs on a fixed set of $n$ vertices, and $ij$-th entry is the number of
ways of removing an edge from $G_j$ and adding an edge to get $G_i$. Note
that the diagonal entry is $m$, since we can add the same edge that is
removed.  Other entries of $B$ are either $0$ or $1$.  The matrix $A$ is
defined similarly for unlabelled graphs with $m$ edges and $n$ vertices.
Thus matrix $A$ is the matrix $\Delta_1$. Matrix
$P$ is defined by indexing the rows by unlabelled graphs and columns by
labelled graphs, and the $ij$-th entry is $1$ if the labelled graph $G_j$
is isomorphic to the unlabelled graph $G_i$.  Other entries are $0$.  As
in [ER], one can verify that $AP=PB$, and every eigenvalue of $A$ is also
an eigenvalue of $B$. But $B=mI+J$, therefore, its eigenvalues are $m +
(m-j)(N-m-j) -j $, where $j \leq $ min$(m,N-m)$.  Thus, eigenvalues of
$(r-1)(2m-r-N)\Delta_0 + \Delta1$ are $ (m-j)(N-m-j+1) +  (r-1)(2m-r-N)$.
If 0 is not an eigenvalue, then $\Delta_{r-1}(X_P-X_Q)=0$, therefore, 
by induction hypothesis, $d_{r-1}(X_P-X_Q) =0 $, and 
$d_r(X_P-X_Q)=0 $ by Kelly's lemma, (see [BH]). For one of the eigenvalues to
be 0, $(r-1)(2m-r-N) \leq 0$. Therefore, $r=1$ (for which the problem is
solved independently in [T]) or $2m \leq N+r$, i.e., $2m-r+1 \leq N+1$.
This contradicts the inequality assumed earlier.

The theorem implies that the $k$-edge deck of a graph can be 
reconstructed from its modified $k$-deck.  

\subsection*{Acknowledgements} This work was done while I was
at Indian Institute of Technology, Guwahati, India.
I would like to thank Philip Maynard for pointing out an 
error in one of the equations in the published version of the paper. 

\section*{References}
\begin{itemize}
\item[{[BH]}] J.A. Bondy  and R. L. Hemminger, Graph reconstruction - a survey,
{\em J. Graph Theory} {\bf 1} no. 3 (1977) 227-268. 
\item[{[BCN]}] A. E. Brouwer, A. M. Cohen and A. Neumaier,
{\em Distance-Regular Graphs}, Springer, Berlin, 1989.
\item[{[ER]}]  M. N. Ellingham and G. F. Royle, Vertex-switching
reconstruction of subgraph numbers and triangle-free graphs, 
{\em J. Combinatorial 
Theory Ser. B}{\bf  54} no. 2 (1992) 167-177.
\item[{[KR]}] I. Krasikov and Y. Roditty, Balance equations for reconstruction
problems, {\em Arch. Math. (Basel)} {\bf 48} no. 5 (1987) 458-464.
\item[{[L]}]  L. Lov\'asz, A note on the line reconstruction problem, {\em 
J. Combinatorial Theory Ser. B} {\bf  13} (1972) 309-310.
\item[{[M]}]  V. M\"uller, The edge reconstruction hypothesis is true for
graphs with more than $n.\log_2 n$ edges, {\em  J. Combinatorial Theory Ser. B}
{\bf 22} no. 3 (1977) 281-283.
\item[{[T]}] B. D. Thatte, A reconstruction problem related to balance
equations-I, {\em Discrete Mathematics} {\bf 176} (1997) 279-284.

\end{itemize}
\end{document}